\documentclass{article}

\usepackage{arxiv}

\usepackage[utf8]{inputenc} 
\usepackage[T1]{fontenc}    
\usepackage{hyperref}       
\usepackage{url}            
\usepackage{booktabs}       
\usepackage{amsfonts}       
\usepackage{nicefrac}       
\usepackage{microtype}      
\usepackage{lipsum}		
\usepackage{graphicx}
\usepackage{doi}

\usepackage{lmodern}
\usepackage[english]{babel}
\usepackage{amsmath,amssymb}
\usepackage{hyperref}
\usepackage{xcolor}
\usepackage{subfig}
\usepackage{graphicx} 
\usepackage{pgf,tikz,svg}

\title{Benchmarking of flatness-based control \\of the heat equation}

\author{\href{https://orcid.org/0000-0003-1762-5659}{\includegraphics[scale=0.06]{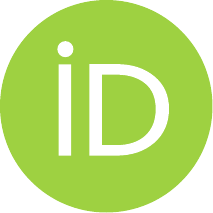}\hspace{1mm} Stephan Scholz}, Lothar Berger \\
	Control and Process Engineering\thanks{Web: \href{https://forschung.rwu.de/forschungsgruppen/control-and-process-engineering}{https://forschung.rwu.de/forschungsgruppen/control-and-process-engineering}},\\
	University of Applied Sciences Ravensburg-Weingarten,\\
	Weingarten, Germany\\
	Email to: \texttt{stephan.scholz@rwu.de} \\
\And
Dirk Lebiedz, \\
Institute of Numerical Mathematics,\\
Ulm University,\\
Ulm, Germany \\
}

\renewcommand{\shorttitle}{Benchmarking of Flatness-based Control}

\hypersetup{
pdftitle={Benchmarking of Flatness-based Control of the Heat Equation},
pdfsubject={},
pdfauthor={Stephan Scholz, Lothar Berger, Dirk Lebiedz},
pdfkeywords={},
}

\begin{document}
\newcommand{\heat}{\vartheta}
\newcommand{\heattx}{\heat(t,x)}
\newcommand{\dotheattx}{\dot{\heat}(t,x)}
\newcommand{\dheattx}{\frac{\partial}{\partial x}  \heat(t,x)}
\newcommand{\ddheattx}{\frac{\partial^2}{\partial x^2}  \heat(t,x)}
\newcommand{\wtx}{w(t,x)}
\newcommand{\dotwtx}{\dot{w}(t,x)}
\newcommand{\opd}{\operatorname{d}}

\definecolor{rwulila}{cmyk}{0.77,0.79,0,0}
\definecolor{rwucyan}{cmyk}{0.88,0,0.11,0}
\definecolor{myblue}{cmyk}{1,0.67,0,0.4}
\definecolor{mywarm}{cmyk}{0,0.29,0.63,0.16}
\definecolor{petrolblue}{cmyk}{0.98,0.06,0.33,0.47}
\definecolor{leafgreen}{cmyk}{0.65,0,1,0}

\definecolor{juliablue}{rgb}{0.251,0.388,0.847}
\definecolor{juliagreen}{rgb}{0.22,0.596,0.149}
\definecolor{juliared}{rgb}{0.796,0.235,0.2}
\definecolor{juliapurple}{rgb}{0.584,0.345,0.698}

\definecolor{plasmacold}{rgb}{0.050383,0.029803,0.527975}
\definecolor{plasmamiddle}{rgb}{0.794549,0.27577,0.473117}
\definecolor{plasmahot}{rgb}{0.940015,0.975158,0.131326}

\pgfdeclarehorizontalshading{plasma}{100bp}{
	color(0bp)=(plasmahot);
	color(25bp)=(plasmahot);
	color(50bp)=(plasmamiddle);
	color(75bp)=(plasmacold);
	color(100bp)=(plasmacold)
}

\maketitle

\begin{abstract}
	Flatness-based control design is a well established method to generate open-loop control signals. Several articles discuss the application of flatness-based control design of (reaction-) diffusion problems for various scenarios. Beside the pure analytical derivation also the numerical computation of the input signal is crucial to yield a reliable trajectory planning. Therefore, we derive the input signal step-by-step and describe the influence of system and controller parameters on the computation of the input signal. In particular, we benchmark the control design of the one-dimensional heat equation with Neumann-type boundary actuation for pure aluminum and steel 38Si7, and discuss the applicability of the found input signals for realistic scenarios.
\end{abstract}


\section{Introduction}

The flatness-based control method is an open-loop technique to steer the system output along a reference trajectory \cite{article:fliess1995flatness}. In case of finite-dimensional linear and nonlinear systems, the input signal $u(t)$ is found by a finite number of derivatives of a (flat) output which equals the reference signal. This approach is extended to infinite-dimensional and distributed parameter systems where theoretically an infinite number of derivatives of output signal $y(t)$ is necessary to compute the input signal \cite{article:laroche2000motion,article:ollivier2001generalization,book:rudolph2003flatness}. However, for practical reasons we can only consider a finite number of derivatives of the output signal. Thus, we need to show that the computation of input signal $u(t)$ converges for a certain number of derivatives of $y(t)$. In general, this estimation of convergence is not trivial because the computation of $u(t)$ depends on system and control parameters. 

In this contribution, we assume a one-dimensional linear heat equation with Neumann-type boundary actuation, as depicted in Figure \ref{fig:one_dim_rod_heat}, and examine the impact of system and control parameters on the computation of input signal $u(t)$. In particular, we compare pure aluminum and steel 38Si7 for this purpose. They differ in their material properties: thermal conductivity $\lambda$, specific heat capacity $c$ and density $\rho$. Regarding the control parameters, we design the reference trajectory as a smooth step function, which is configured by the transition time and the steepness \cite{article:utz2010trajectory}. 

In Section \ref{sec:flatness_based_control}, we introduce the flatness-based modeling for the one-dimensional heat equation and derive input signal $u(t)$. The influence of the system parameters are analyzed in Section \ref{sec:influence_system}. The trajectory planning problem and the subsequent discussion of the control parameters are described in Section \ref{sec:trajectory_planning} and \ref{sec:influence_control}, respectively. Finally, we present the simulation results of the open-loop control system and review the applicability for realistic scenarios in Section \ref{sec:simulation_results}.

\section{Flatness-based Control}
\label{sec:flatness_based_control}

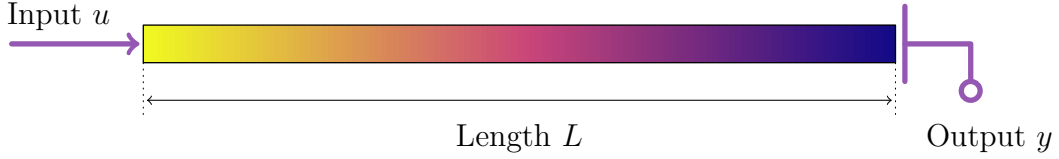
\begin{figure}[t]
	\centering
	\begin{tikzpicture}[scale=2.5,line cap=round,line join=round,x=1.0cm,y=1.0cm]
		\shade[shading=plasma] ((0,0) rectangle (4,0.2);
		\draw [line width=0.5, color=black] (0,0) rectangle (4,0.2);
		\draw [->, line width=2., color=juliapurple] (-0.7, 0.1) -- (-0.02, 0.1);
		\draw[color=black] (-0.45, 0.25) node {\large Input $u$};
		\draw [line width=2., color=juliapurple] (4.05,0.3) -- (4.05,-0.1);
		\draw [line width=2., color=juliapurple] (4.05,0.1) -- (4.4, 0.1) -- (4.4, -0.1);
		\draw [line width=2., color=juliapurple] (4.4,-0.15) circle [radius=0.05];
		\draw[color=black] (4.5, -0.4) node {\large Output $y$};
		\draw [dotted] (0,0) -- (0.0, -0.3);
		\draw [dotted] (4,0) -- (4.0, -0.3);
		\draw [<->] (0.02, -0.2) -- (3.98, -0.2);
		\draw[color=black] (2., -0.4) node {\large Length $L$};
	\end{tikzpicture}
	\caption{One-dimensional rod with heat input (left) and temperature measurement (right).}
	\label{fig:one_dim_rod_heat}
\end{figure}

We assume a one-dimensional heat equation with Neumann-type boundary actuation on the left side and a thermally insulated right side as
\begin{align}
	\dotheattx =&~ \alpha \ddheattx \quad \text{,} \quad (t,x) \in (0,T) \times (0,L) \text{,}  \label{eq:fbc_heat_eq} \\
	u(t) =&~ \lambda \left. \dheattx \cdot \vec{n}_{0} \right\rvert_{x=0}   \text{,} \label{eq:fbc_heat_eq_bc_left} \\
	0 =&~ \lambda \left. \dheattx \cdot \vec{n}_{L} \right\rvert_{x=L}   \label{eq:fbc_heat_eq_bc_right}
\end{align}
where the outer normal vectors are known as $\vec{n}_{0}=-1$ and $\vec{n}_{L}=1$. The heat conduction model is portrayed in Figure \ref{fig:one_dim_rod_heat}. Here, we denote the temperature  as $\heat$, the thermal conductivity as $\lambda$ and the diffusivity as $\alpha=\frac{\lambda}{c~\rho}$ with specific heat capacity $c$ and density $\rho$. The boundary conditions (\ref{eq:fbc_heat_eq_bc_left},~\ref{eq:fbc_heat_eq_bc_right}) are rough simplifications because realistic scenarios may treat linear heat transfer and nonlinear heat radiation, see \cite[chapter 1.1]{book:baehr2013heat} and \cite[chapter 1.3]{book:lienhard2020heat}. Such realistic boundary conditions are described in \cite{article:herzog2023an,article:scholz2020modeling,article:scholz2022hestia}. They are not treated here because they lead to much more complex mathematical discussions which are out of scope of this contribution. The initial temperature distribution is defined by 
\begin{align*}
	\heat(0,x) =&~ \heat_{0}(x) \quad \text{for} \quad x \in [0,L]
\end{align*}
and the temperature is measured on the right boundary side as
\begin{align}
	y(t) = \heat(t,L) \text{.} \label{eq:fbc_heat_eq_output}
\end{align}

As known from the literature \cite{article:laroche2000motion,article:ollivier2001generalization,book:rudolph2003flatness}, the heat equation can be represented by a power series approach. So, we define the power series
\begin{align*}
	w(t,x) := \sum_{i=0}^{\infty} w_{i}(t) \frac{(L-x)^{i}}{i!}
\end{align*}
and find its derivatives with respect to position $x$ as
\begin{align}
	\frac{\partial}{\partial x} w(t,x) =&~ -\sum_{i=0}^{\infty} w_{i+1}(t) \frac{(L-x)^{i}}{i!} \quad \text{and} \label{eq:fbc_power_series_dx} \\
	\frac{\partial^2}{\partial x^2} w(t,x) =&~ \sum_{i=0}^{\infty} w_{i+2}(t) \frac{(L-x)^{i}}{i!} \text{.} \nonumber 
\end{align}
We model heat equation \eqref{eq:fbc_heat_eq} in terms of
\begin{align*}
	\dot{w}(t,x) =&~ \alpha ~ \frac{\partial^2}{\partial x^2} w(t,x) \text{,}
\end{align*}
identify both sides by its power series expressions as
\begin{align*}
	\sum_{i=0}^{\infty} \dot{w}_{i}(t) \frac{(L-x)^{i}}{i!} ~=~\alpha~ \sum_{i=0}^{\infty} w_{i+2}(t) \frac{(L-x)^{i}}{i!} 
\end{align*}
and yield identity
\begin{align}
	\dot{w}_{i}(t) = \alpha~ w_{i+2}(t) \text{.} \label{eq:fbc:identity_seq}
\end{align}
Next, we apply the information of both boundary sides on identity \eqref{eq:fbc:identity_seq} to derive the input signal. Firstly, we consider the output signal \eqref{eq:fbc_heat_eq_output} as
\begin{align*}
	y(t) = w(t,L) = \sum_{i=0}^{\infty} w_{i}(t) \frac{0^{i}}{i!} = w_{0}(t)
\end{align*}
which implies $\frac{d^{i}}{d t^{i}} y(t) = \frac{d^{i}}{d t^{i}} w_{0}(t) = \alpha^{i} w_{2i}$ with identity \eqref{eq:fbc:identity_seq}. Secondly, the boundary condition on the right side \eqref{eq:fbc_heat_eq_bc_right} is formulated as
\begin{align*}
	\lambda \frac{\partial}{\partial x} w(t,L) = -\lambda \sum_{i=0}^{\infty} w_{i+1}(t) \frac{0^{i}}{i!} = - \lambda w_{1}(t) = 0
\end{align*}
and we find $ \frac{d^{i}}{d t^{i}} w_{1}(t) = \alpha^{i} w_{2i+1} \equiv 0$. Thus, identity \eqref{eq:fbc:identity_seq} is split up into the sequences
\begin{align*}
	w_{2i}(t) =~ \alpha^{-i} y^{(i)}(t) \quad \text{and} \quad
	w_{2i+1}(t) =~ 0 \quad \text{for all} ~ i\in\{0,1,\ldots,\infty\}\text{.} 
\end{align*}
We insert Eq. \eqref{eq:fbc_power_series_dx} in the definition of boundary actuation \eqref{eq:fbc_heat_eq_bc_left} to derive the input signal $u(t)$ as
\begin{align}
	u(t) =&~- \lambda \frac{\partial}{\partial x} w(t,0) = \lambda \sum_{i=0}^{\infty} w_{i+1}(t) \frac{L^{i}}{i!} \text{.} \label{eq:fbc_input_signal_prep}
\end{align}
We know that only $w_{2i}$ entries are not zero. So, we map index  $i \rightarrow 2i+1$ in Eq. \eqref{eq:fbc_input_signal_prep} and find the input signal 
\begin{align}
	u(t) =~ \lambda \sum_{i=0}^{\infty} \frac{L^{2i+1}}{\alpha^{i+1}} \frac{1}{(2i+1)!} ~ y^{(i+1)}(t) \text{.} \label{eq:fbc_input_signal}
\end{align}
Input signal \eqref{eq:fbc_input_signal} consists of the derivatives of the output signal $y$ and the system parameters: length $L$, thermal conductivity $\lambda$ and diffusivity $\alpha$. In the design of control algorithms, we are usually only able to specify the desired output signal and its derivatives, as explained in Section \ref{sec:trajectory_planning}. The system parameters are usually predefined (and are not modifiable) by the choice of material and geometry. Nonetheless, the system parameters have a significant influence on the computation of input signal \eqref{eq:fbc_input_signal}, as described below.

\section{Influence of System Parameters}
\label{sec:influence_system}

We are interested in the sequence values of series \eqref{eq:fbc_input_signal} because for implementation reasons we need to know how much memory has to be reserved for the computation of $u$ and at which iteration $i$ the summation can be stopped. We extract the sequence 
\begin{align}
	\eta_{i} =  \frac{L^{2i+1}}{\alpha^{i+1}} \frac{1}{(2i+1)!} \label{eq:fbc_eta_sequence}
\end{align}
from the power series in Eq. \eqref{eq:fbc_input_signal}, and analyze the influence of length $L$ and diffusivity $\alpha$ on sequence $\eta_{i}$ in this section. The influence of the desired output signal $y(t)$ and its derivatives is examined in Section \ref{sec:influence_control}. 

Sequence $\eta_{i}$ is positive for all $i\in\{0,1,\ldots,\infty\}$ as we assume $L>0$, $\alpha>0$, and it scales the derivatives $y^{(i+1)}$. We rescale the sequence \eqref{eq:fbc_eta_sequence} as 
\begin{align*}
	\tilde{\eta}_{i} := \left(\frac{L^{2}}{\alpha}\right)^{i+1} \frac{1}{(2i+1)!} =  \frac{\gamma^{i+1}}{(2i+1)!} = L ~ \eta_{i} 
\end{align*}
with $\gamma := \frac{L^2}{\alpha}$ to show that $\eta_{i}$ and $\tilde{\eta}_{i}$ increase up to some index $i$ and decrease afterwards to zero. Coefficient $\gamma$ is also known as part of the Fourier number $\operatorname{Fo}=\frac{\alpha}{L^2}~t= \frac{t}{\gamma}$ with time $t\geq0$, see also \cite[page 116]{book:baehr2013heat} and \cite[page 195]{book:lienhard2020heat}. When we increase iterator $i$ by one, then we yield
\begin{align*}
	\tilde{\eta}_{i+1} = \frac{\gamma^{[i+1]+1}}{(2[i+1]+1)!} = \frac{\gamma^{i+1}}{(2i+1)!} \frac{\gamma}{(2i+2)(2i+3)} = \tilde{\eta}_{i} ~  \beta_{i}
\end{align*}
where $\beta_{i}=\frac{\gamma}{(2i+2)(2i+3)}$. Consequently, we notice 
\begin{align*}
	\frac{\tilde{\eta}_{i+1}}{\tilde{\eta}_{i}} >~ 1 \quad \Leftrightarrow \quad  \beta_{i}> 1 \quad \text{and} \quad	
	\frac{\tilde{\eta}_{i+1}}{\tilde{\eta}_{i}} <~ 1 \quad \Leftrightarrow \quad \beta_{i}< 1 \text{.}
\end{align*}
Due to the definition of $\tilde{\eta}$, this concept holds also for the original sequence \eqref{eq:fbc_eta_sequence} as $\eta_{i+1} = \beta_{i} ~ \eta_{i}$. Thus, the maximum value of $\tilde{\eta}_{i}$ and $\eta_{i}$ and its corresponding iterations $i_{max}$ depend only on $\gamma$. For example, if we assume $\gamma=100$ then $\gamma < (2i+2)(2i+3)$ holds for $i\in\{1,2,3\}$ and we find the maximum value $\tilde{\eta}_{4} = \frac{100^{5}}{9!}\approx 27557$. 

\subsection*{Example: Comparison of Aluminum and Steel 38Si7}

\def\arraystretch{1.3}
\begin{table}[b]
	\begin{center}
		\caption{\MakeUppercase{Physical Properties}}
		\begin{tabular}{|l | l l l l|}
			\hline
			& $\lambda$ & $\rho$ & $c$ & $\alpha=\frac{\lambda}{\rho c}$ \\
			Aluminum & $237$ & $2700$ & $900$ & $9.75\cdot10^{-5}$\\
			Steel 38Si7 & $40$ & $7800$ & $460$ & $1.11\cdot10^{-5}$\\
			\hline
		\end{tabular}
		\label{table:physical_prop_al_steel}
	\end{center}
\end{table}

\begin{figure}[!t]
	\centering
	\subfloat[Sequence $\eta_{i}$]{\includegraphics[width=0.48\columnwidth]{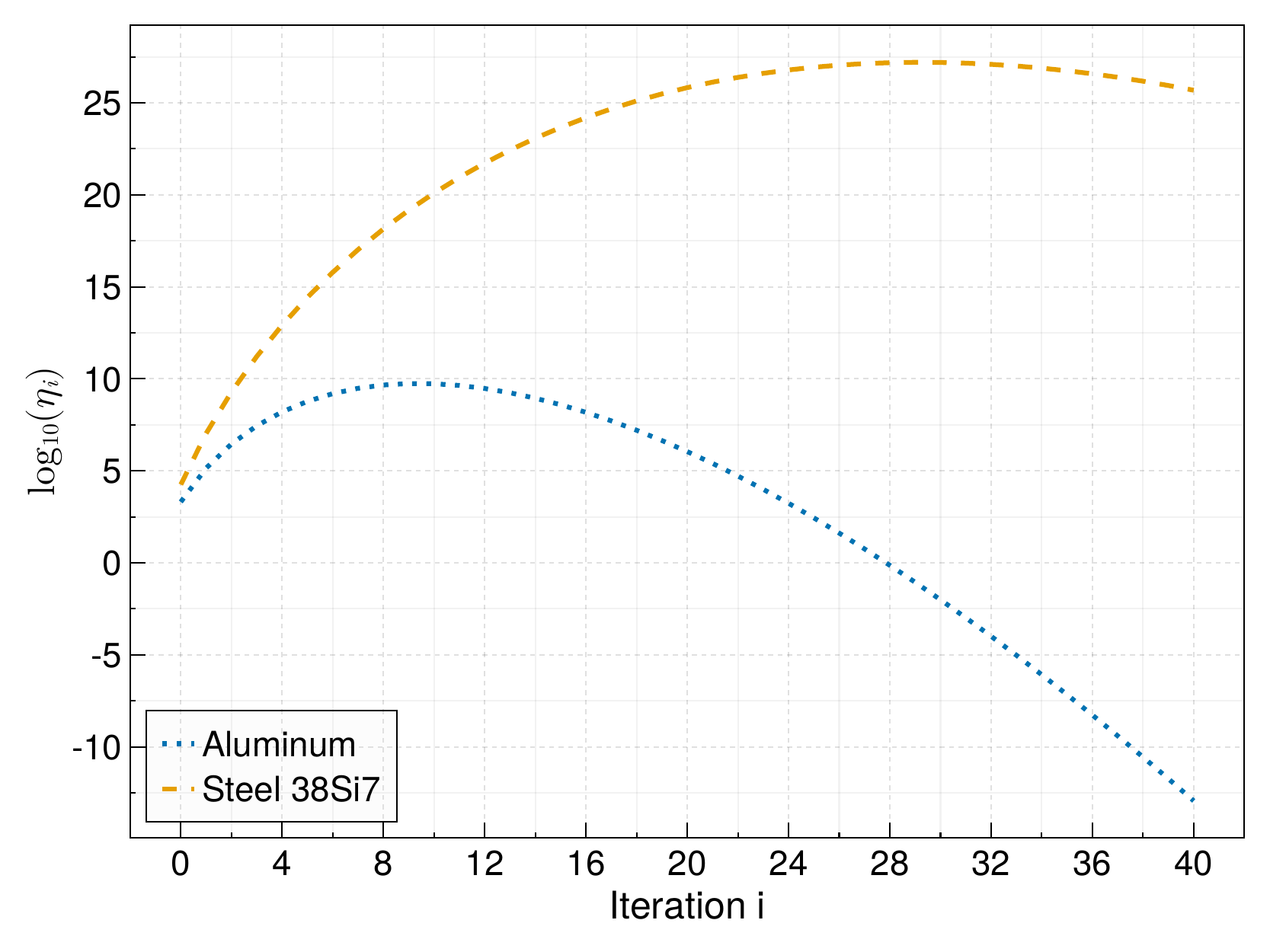}}
	\subfloat[Ratio  $\eta_{i+1}/\eta_{i}$]{\includegraphics[width=0.48\columnwidth]{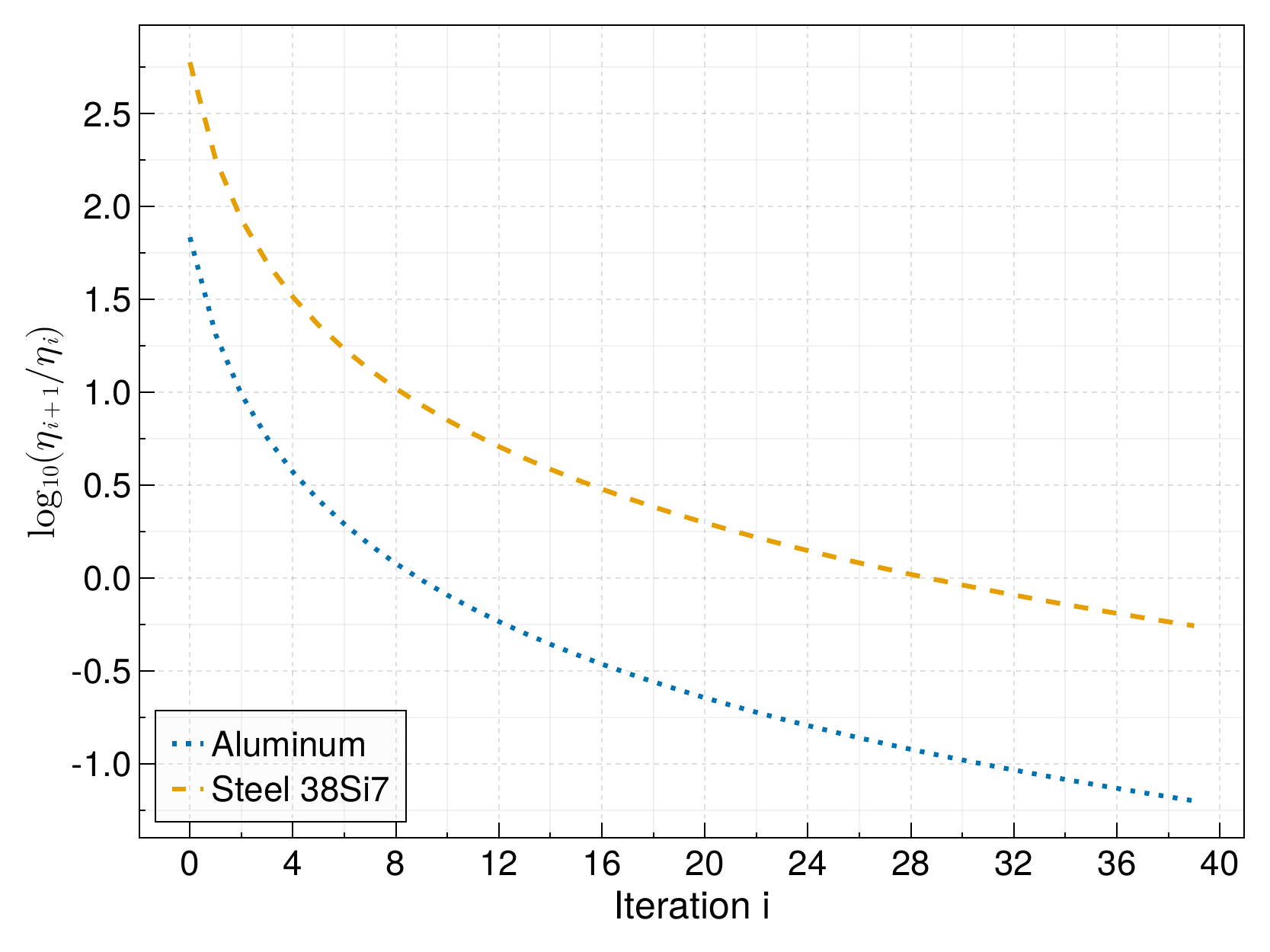}}
	\caption{Sequence $\eta_{i}$ (left) and ratio $\frac{\eta_{i+1}}{\eta_{i}}$ (right) for aluminum and steel 38Si7 in semi-logarithmic scaling.}
	\label{fig:sequence_eta_al_steel}
\end{figure}

For our numerical evaluations we consider a rod of length $L=0.2$ for two case scenarios: a rod made of pure aluminum \cite{online:periodic2023aluminium} and a rod made of steel 38Si7 \cite{online:ovako202338Si7}. The physical properties of both materials are listed in \mbox{Table \ref{table:physical_prop_al_steel}}. For aluminum we have $\gamma_{al} \approx 410$ and for steel 38Si7 we have $\gamma_{st} \approx 3588$.
We portray in \mbox{Figure \ref{fig:sequence_eta_al_steel}} in semi-logarithmic scaling the sequences $\eta_{al,i}$ (aluminum) and $\eta_{st,i}$ (steel 38Si7), and their ratios $\frac{\eta_{al,i+1}}{\eta_{al,i}}$ and $\frac{\eta_{st,i+1}}{\eta_{st,i}}$. These ratios describe the evolution of the sequences for an increasing iteration $i$. We find that inequality $\frac{\eta_{i+1}}{\eta_{i}} >~ 1$ or equally $\log_{10}\left(\frac{\eta_{i+1}}{\eta_{i}}\right) >~ 0$ holds in case of aluminum for $i\in\{1,\ldots,8\}$ and in case of steel $i\in\{1,\ldots,28\}$. Thus, the maximum values of $\eta_{i}$ for aluminum and steel are calculated by
\begin{align*}
	\eta_{al,9}=~\frac{L^{19}}{\alpha_{al}^{10} ~ 19!} \approx 5.53 \cdot 10^{9} \quad \text{and}  \quad
	\eta_{st,29}=~\frac{L^{59}}{\alpha_{st}^{30} ~ 59!} \approx 1.59 \cdot 10^{27} \text{.}
\end{align*}
As both sequences $\eta_{al,i}$ and $\eta_{st,i}$ reach such enormous maximum values, computational issues related to big numbers and data types have to be considered in the implementation process.

Moreover, sequence $\log_{10}(\eta_{al,i})$ drops below zero for $i>27$: $\eta_{al,28}\approx0.73$,  $\log_{10}(\eta_{al,28})\approx -0.13$; and $\log_{10}(\eta_{st,i})$ drops below zero for $i>82$: $\eta_{st,83}\approx0.13$, $\log_{10}(\eta_{st,83})\approx -0.87$ (not displayed in Figure \ref{fig:sequence_eta_al_steel}). 

\section{Trajectory Planning}
\label{sec:trajectory_planning}
As introduced in Section \ref{sec:flatness_based_control}, we assume a temperature measurement $y$ on the right boundary side of the rod at $x=L$, see Eq. \eqref{eq:fbc_heat_eq_output}. This measured temperature shall be steered along a predefined trajectory. According to \cite{book:rudolph2003flatness,article:utz2010trajectory}, we consider a transition from the initial fixed operating point $y_{0}$ to the next one $y_{f}$ as
\begin{align}
	y(t) = y_{0} + \Delta y ~ \Phi_{\omega,T}(t)
\end{align}
where $\Delta y = y_{f} - y_{0}$, and with transition function
\begin{align}
	\Phi_{\omega,T}(t) = 
	\begin{cases}
		0 & \text{for}~ t \leq 0 \text{,} \\
		1 & \text{for}~ t \geq T \text{,} \\
		\frac{\int_{0}^{t} \Omega_{\omega,T}(\tau) \opd\tau}{\hat{\Omega}_{\omega,T}} & \text{for}~  t \in (0,T) \text{.}
	\end{cases} \label{eq:fbc_transition}
\end{align}
The transition or smooth step function \eqref{eq:fbc_transition} consists of the bump function
\begin{align}
	\Omega_{\omega,T}(t) = 
	\begin{cases}
		0 & \text{for}~ t \notin [0,T] \text{,}\\
		\exp\left(-1/\left([1 - \frac{t}{T}] \frac{t}{T} \right)^\omega \right) & \text{for}~ t \in (0,T) 
	\end{cases}  \label{eq:fbc_internal_fun}
\end{align}
and its integral
\begin{align}
	\hat{\Omega}_{\omega,T} :=&~ \int_{0}^{T} \Omega_{\omega,T}(\tau) \opd\tau \text{.}  \label{eq:omega_integral}
\end{align}

\begin{figure}[!t]
	\centering
	\subfloat[Transition function]{\includegraphics[width=0.48\columnwidth]{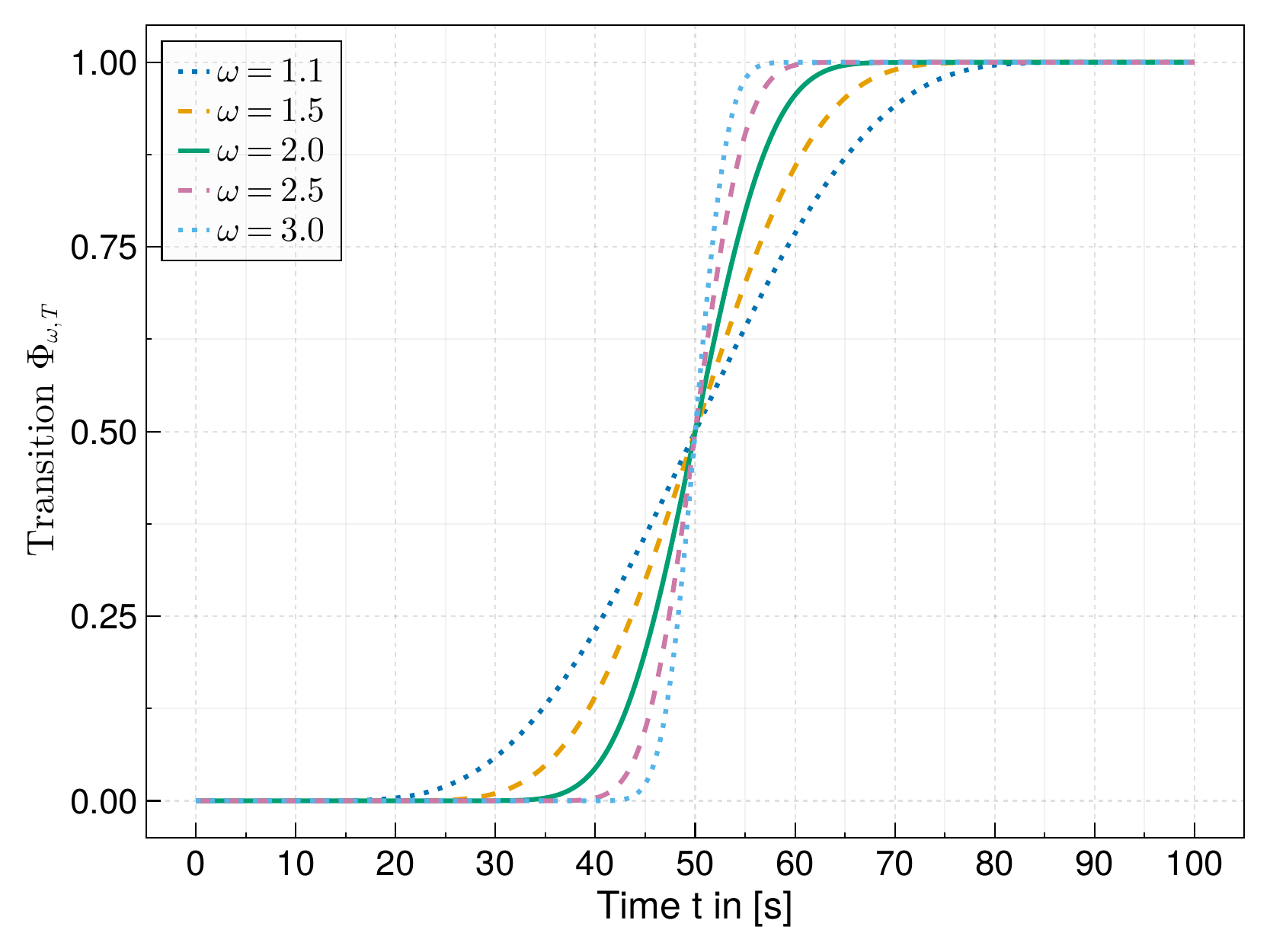}}
	\subfloat[First  derivative]{\includegraphics[width=0.48\columnwidth]{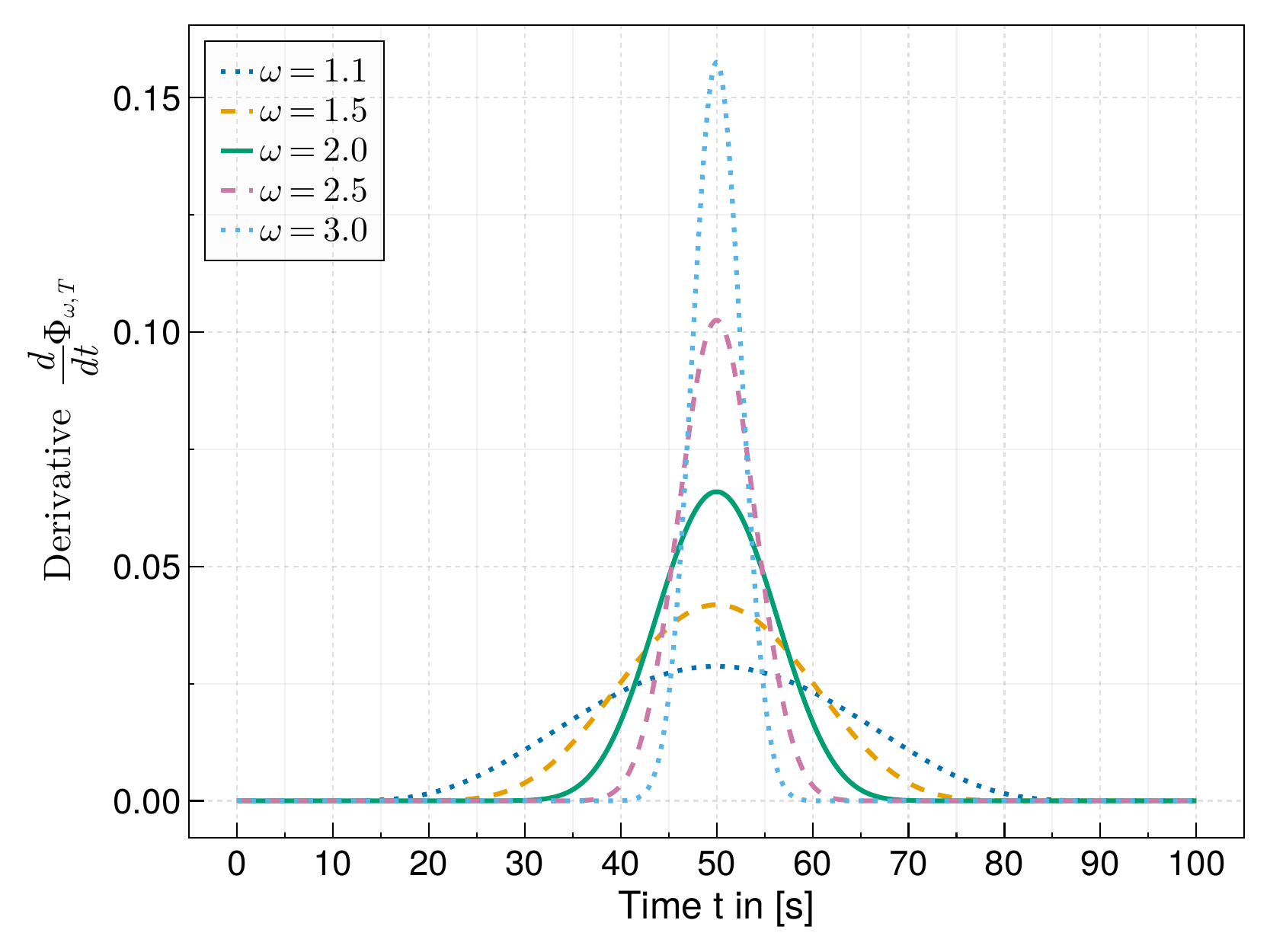}}
	\caption{Transition $\Phi_{\omega,T}$ as in Eq. \eqref{eq:fbc_transition} (left)  and its first derivative (right) with $T=100$ seconds and varying $\omega \in \{1.1, 1.5, 2.0, 2.5, 3.0\}$.}
	\label{fig:trajectory_vary_w}
\end{figure}

Parameter $\omega$ steers the steepness of transition $\Phi_{\omega,T}$ and is chosen such that the Gevrey order $go=1 + \frac{1}{\omega}<2$ or equally $\omega > 1$. More information about transition \eqref{eq:fbc_transition} and bump function \eqref{eq:fbc_internal_fun} can be found in \cite{article:laroche2000motion}. A small value of $\omega$, e.g. $\omega=1.1$ means a rather flat transition, whereas a large value, e.g. $\omega=3.0$ means a quite steep transition, as depicted in Figure \ref{fig:trajectory_vary_w}. We need to find the derivatives 
\begin{align}
	\frac{d^{i}}{dt^{i}} y(t) = \Delta y ~  \Phi_{\omega,T}^{(i)}(t) \label{eq:y_der_phi}
\end{align}
for the computation of input signal $u(t)$ in Eq. \eqref{eq:fbc_input_signal}. From Eq. \eqref{eq:fbc_transition}, we derive the derivatives of the transition function as 
\begin{align}
	\Phi_{\omega,T}^{(i)}(t) =&~
	\begin{cases}
		\frac{\Omega_{\omega,T}^{(i-1)}(t)}{\hat{\Omega}_{\omega,T}} & \text{for} \quad t \in (0,T) \text{,} \\
		0 & \text{for} \quad t \notin (0,1) \text{.}
	\end{cases}
	  \label{eq:phi_der_omega}
\end{align}
Trajectory $\Phi_{\omega,T}(t)$ and its first derivative are portrayed in Figure \ref{fig:trajectory_vary_w} for varying $\omega \in \{1.1, 1.5, 2.0, 2.5, 3.0\}$. In this contribution, we compute the derivatives $\Omega_{\omega,T}^{(i)}$ with the \textsc{Julia} library \textit{BellBruno.jl} \cite{software:scholz2023bellbruno}. Alternatively, the derivatives may be computed with computer algebra systems, see e.g. this \textsc{Matlab} implementation \cite{software:fischer2021coni}.

\section{Influence of Control Parameters}
\label{sec:influence_control}
\begin{figure}[!t]
	\centering	
	\subfloat[Fixed $\omega = 2.0$, varying $T$]{	\includegraphics[width=0.48\columnwidth]{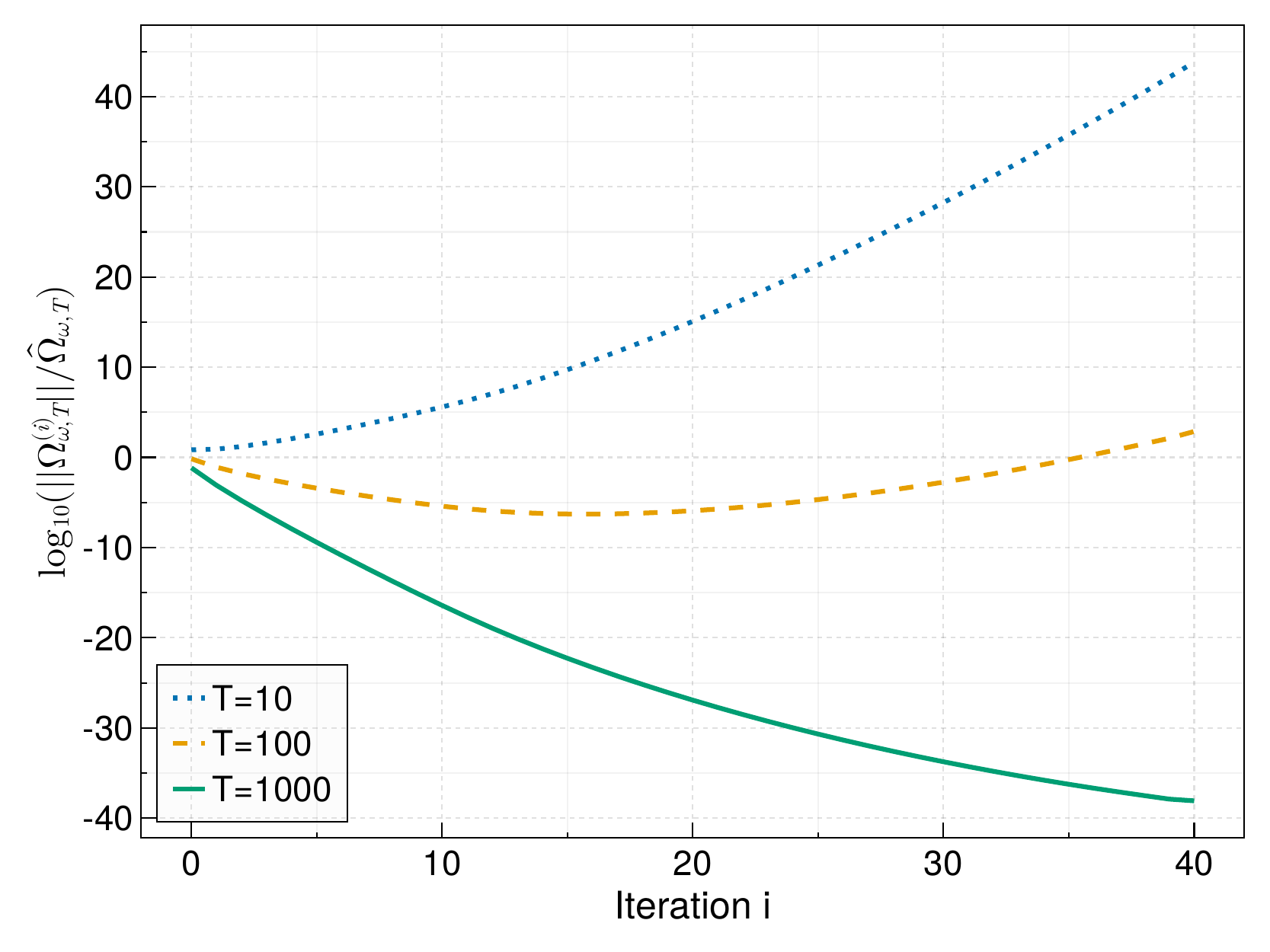}}
	\subfloat[Fixed $T = 1000$, varying $\omega$]{\includegraphics[width=0.48\columnwidth]{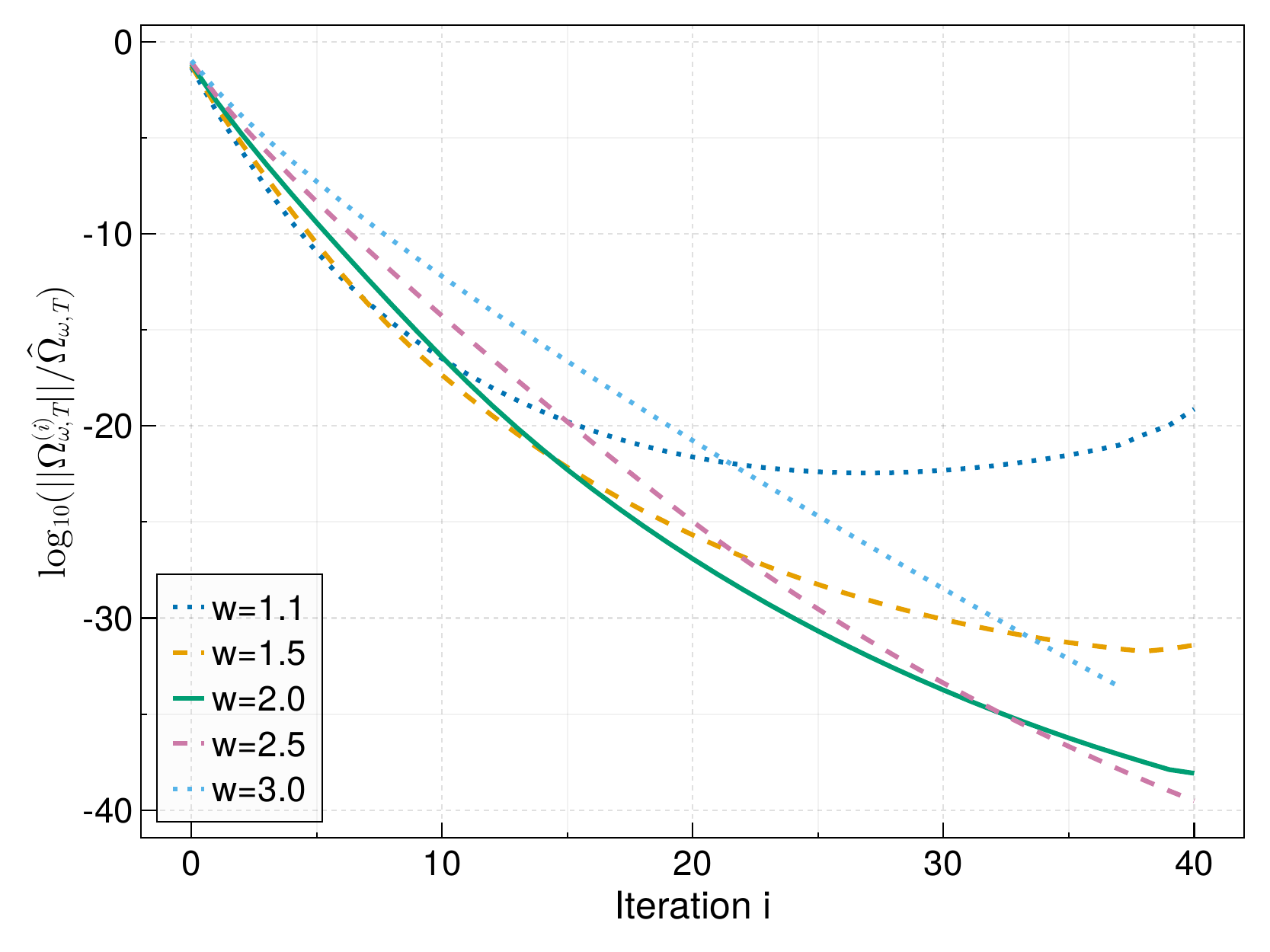}}
	\caption{Normed derivatives $\frac{\left\lVert\Omega_{\omega,T}^{(i)}(t) \right\rVert_{L^2}}{\hat{\Omega}_{\omega,T}}$ in semi-logarithmic scaling for fixed $\omega = 2.0$ (left), and fixed $T = 1000$ (right).}
	\label{fig:norm_omega}
\end{figure}

The configuration of transition $\Phi_{\omega,T}$ and its derivatives are mainly driven by two parameters: final time $T$ and exponent $\omega$. In this section, we apply the $L^{2}$ norm 
\begin{align*}
	\lVert f \rVert_{L^2} = \sqrt{\int_{0}^{T} \lvert f(t) \rvert^2 dt}
\end{align*}
on $\frac{d^{i}}{dt^{i}} \Omega_{\omega,T}(t)$ to unveil the influence of final time $T$ and exponent $\omega$ on the computation of input signal \eqref{eq:fbc_input_signal}. We note the input signal with sequence $\eta_{n}$ as 
\begin{align*}
	u(t) =&~ \lambda \sum_{i=0}^{\infty} \eta_{i} ~ y^{(i+1)}(t) 
	~=~ \frac{\lambda ~ \Delta y}{	\hat{\Omega}_{\omega,T}} ~ \sum_{i=0}^{\infty} \eta_{i} ~ \Omega_{\omega,T}^{(i)}(t)
\end{align*}
and consider the identities (\ref{eq:y_der_phi},\ref{eq:phi_der_omega}) to find the $L^2$ norm of $u(t)$ as
\begin{align*}
	\lVert u(t) \rVert_{L^2} =&~ \left\lVert  \frac{\lambda ~ \Delta y}{\hat{\Omega}_{\omega,T}} ~ \sum_{i=0}^{\infty} \eta_{i} ~ \Omega_{\omega,T}^{(i)}(t) \right\rVert 
	~\leq~ \left\lvert \Delta y  \right\rvert ~ \frac{\lambda}{\hat{\Omega}_{\omega,T}} ~ \sum_{i=0}^{\infty} \eta_{i} ~ \left\lVert  \Omega_{\omega,T}^{(i)}(t) \right\rVert
\end{align*}
where we assume $\lambda > 0,~ \hat{\Omega}_{\omega,T} > 0,~ \eta_{i} > 0$.  We see that the power series is mainly driven by $\eta_{i}$ (as discussed in Section \ref{sec:influence_system}) and derivatives $ \Omega_{\omega,T}^{(i)}(t)$. Therefore, we are able to describe the quantitative behavior of the input signal by evaluating sequence
\begin{align}
	\mu_{i} := \lambda ~ \lvert \Delta y \rvert ~ \eta_{i} ~ \frac{\left\lVert\Omega_{\omega,T}^{(i)}(t) \right\rVert_{L^2}}{\hat{\Omega}_{\omega,T}}\text{.} \label{eq:seq_mu}
\end{align} 
We remark that sequence \eqref{eq:seq_mu} depends mainly on $\eta_{i}$ and on $\frac{\left\lVert\Omega_{\omega,T}^{(i)}(t) \right\rVert_{L^2}}{\hat{\Omega}_{\omega,T}}$ which represent the system and control parameters, respectively. The normed derivatives $\frac{\left\lVert\Omega_{\omega,T}^{(i)}(t) \right\rVert_{L^2}}{\hat{\Omega}_{\omega,T}}$ are portrayed in Figure \ref{fig:norm_omega} in semi-logarithmic scaling for two scenarios: fixed $\omega=2$ and varying $T\in\{10, 100, 1000\}$; and fixed $T = 1000$ and varying $\omega \in \{1.1, 1.5, 2.0, 2.5, 3.0\}$. Figure \ref{fig:norm_omega} illustrates, that an increasing value of final time $T$ leads to a significant reduction of $\frac{\left\lVert\Omega_{\omega,T}^{(i)}(t) \right\rVert_{L^2}}{\hat{\Omega}_{\omega,T}}$ per iteration $i$. In contrast to this situation, the effect of a variation of steepness $\omega$ may not be so clear here.

\begin{figure}[!t]
	\centering
	\subfloat[Sequence $\mu_{i}$ from Eq. \eqref{eq:seq_mu}]{\includegraphics[width=0.48\columnwidth]{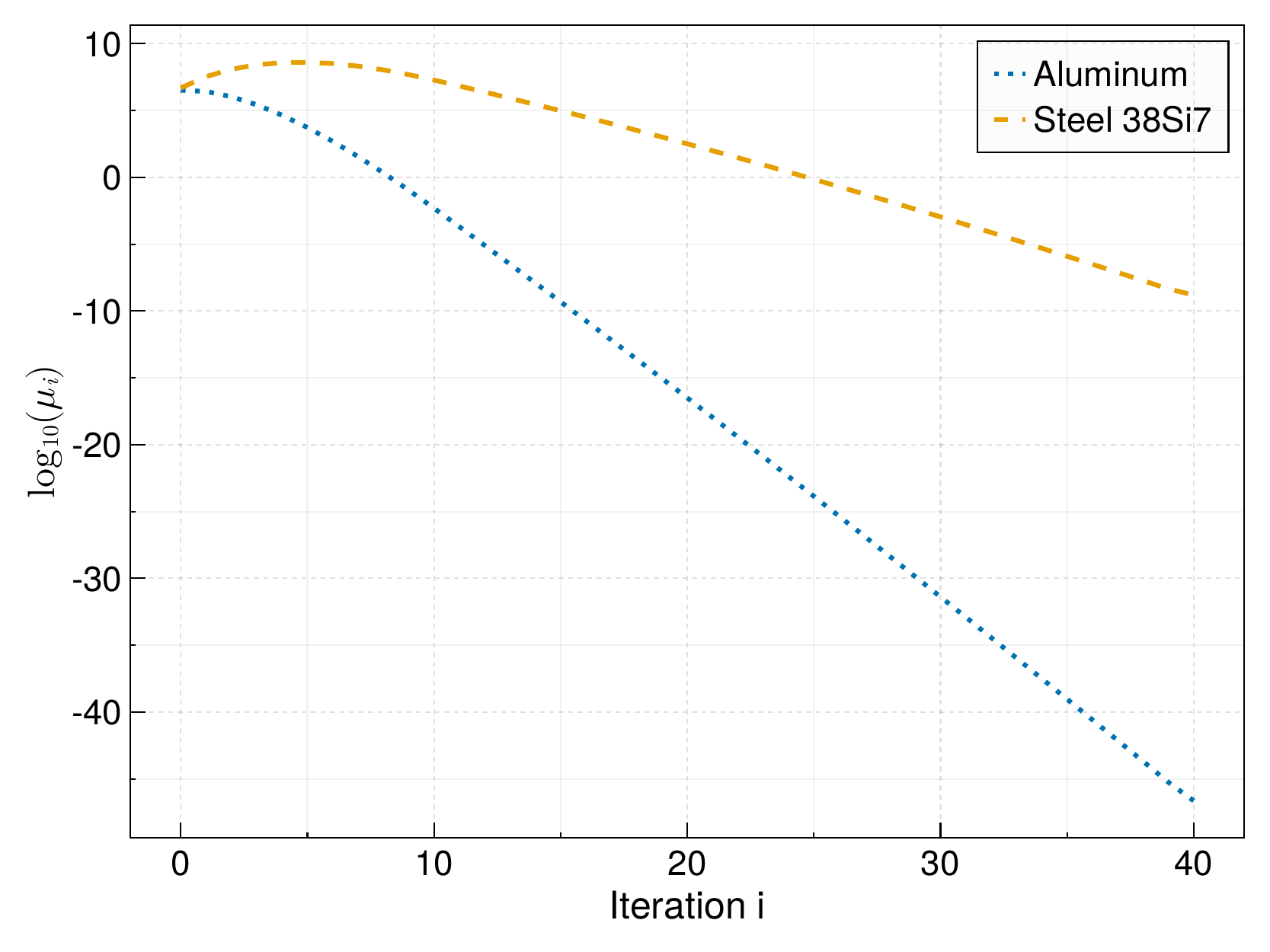}}
	\subfloat[Ratio $\hat{r}_{i}$ from Eq. \eqref{eq:ratio_r}]{\includegraphics[width=0.48\columnwidth]{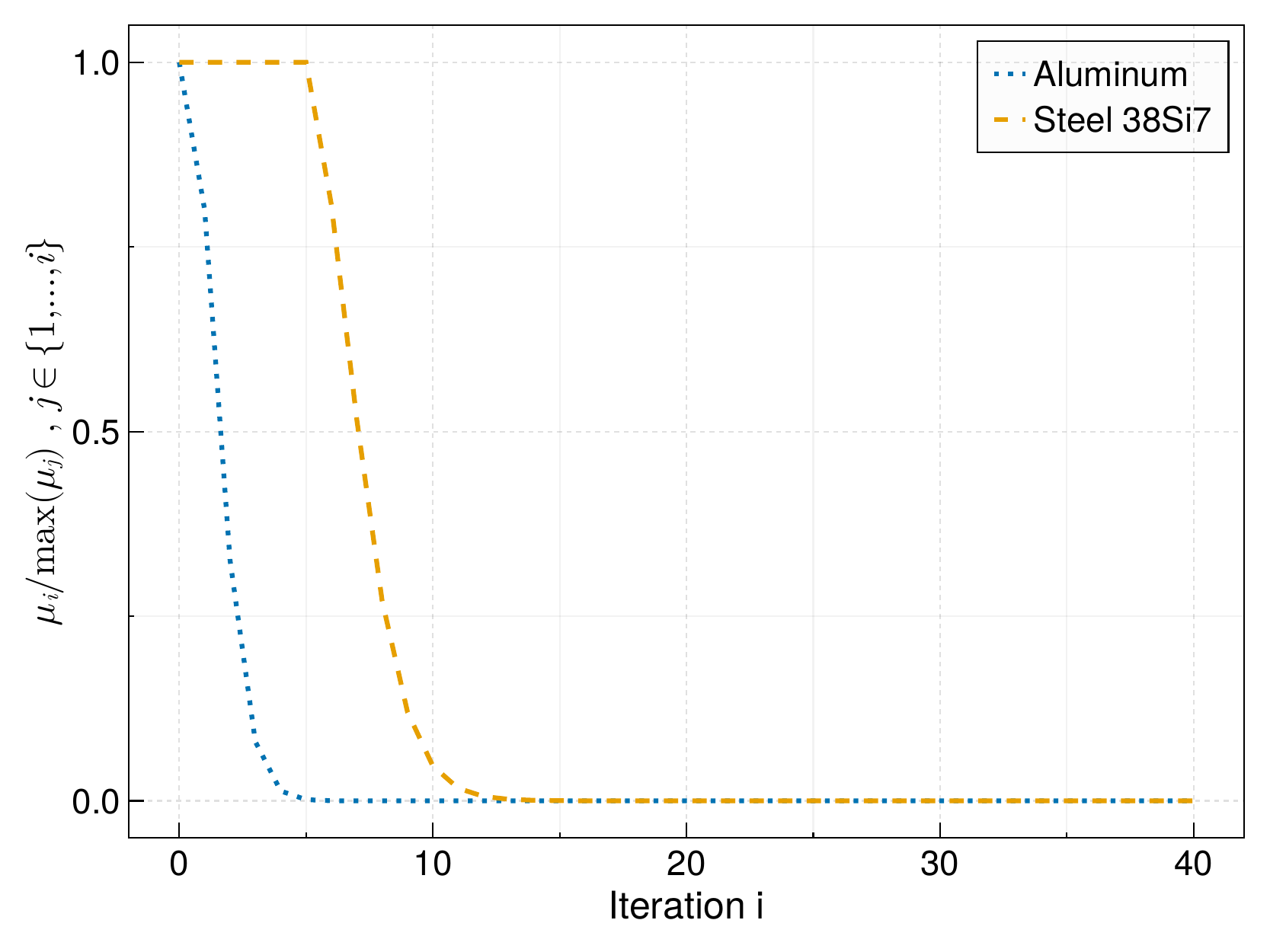}}
	\caption{Sequence $\mu_{i}$ in semi-log. scaling and ratio $\hat{r}_{i}$ for $\omega=2.0$ and $T=1000$.}
	\label{fig:norm_mu}
\end{figure}

Furthermore, we take advantage of sequence $\mu_{i}$ to find a suitable maximum iteration number $i_{max}$ to terminate the power series in Eq. \eqref{eq:fbc_input_signal}. The sequence $\mu_{i}$ is computed in one case for aluminum with $\eta_{al,i}$ and in the other case for steel 38Si7 with $\eta_{st,i}$, as introduced in Section \ref{sec:influence_system}. The different values of $\eta_{i}$ for aluminum and steel 38Si7, as depicted in Figure \ref{fig:sequence_eta_al_steel}, lead to different values of $\mu_{i}$. The sequence $\mu_{i}$ approaches zero \textit{faster} in case of aluminum than in case of steel 38Si7 as depicted in Figure \ref{fig:norm_mu} (a). 
We introduce the ratio 
\begin{align}
	\hat{r}_{i} := \frac{\mu_{i}}{\max\limits_{j \in \{1,\ldots,i\}} \mu_{j}} \label{eq:ratio_r}
\end{align}
to compare the value of $\mu_{i}$ with its predecessors, and we portray $\hat{r}_{i}$ in Figure \ref{fig:norm_mu} (b). The ratio $\hat{r}_{i}$ approaches zero in case of aluminum for $i>5$ and in case of steel 38Si7 for $i>12$. The data of sequence $\mu_{i}$ and ratio $\hat{r}_{i}$ in Figure \ref{fig:norm_mu} unveil that the series in input signal \eqref{eq:fbc_input_signal} needs more iterations for steel 38Si7  to converge compared to aluminum. More iterations entail higher derivatives $\Omega_{\omega,T}^{(i)}(t)$, which tend to more oscillatory behavior. The series elements are also larger for steel 38Si7 and hence we expect higher amplitudes of the input signal. Therefore, the generation of the input signal for steel 38Si7 may need notable more energy and data memory than the input signal for aluminum. 

\begin{figure}[!t]
	\centering
	\subfloat[Aluminum]{\includegraphics[width=0.48\columnwidth]{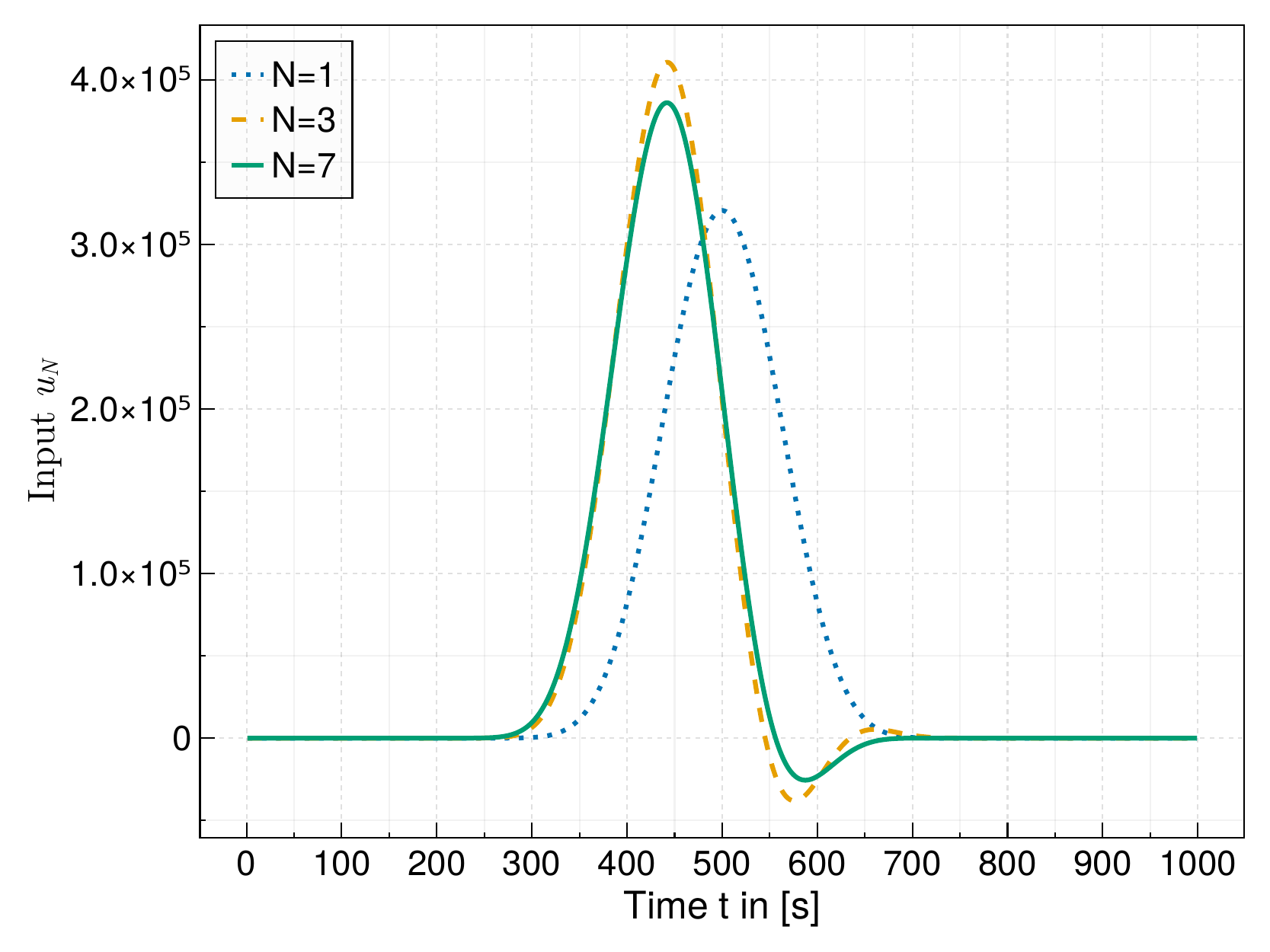}}
	\subfloat[Steel 38Si7]{\includegraphics[width=0.48\columnwidth]{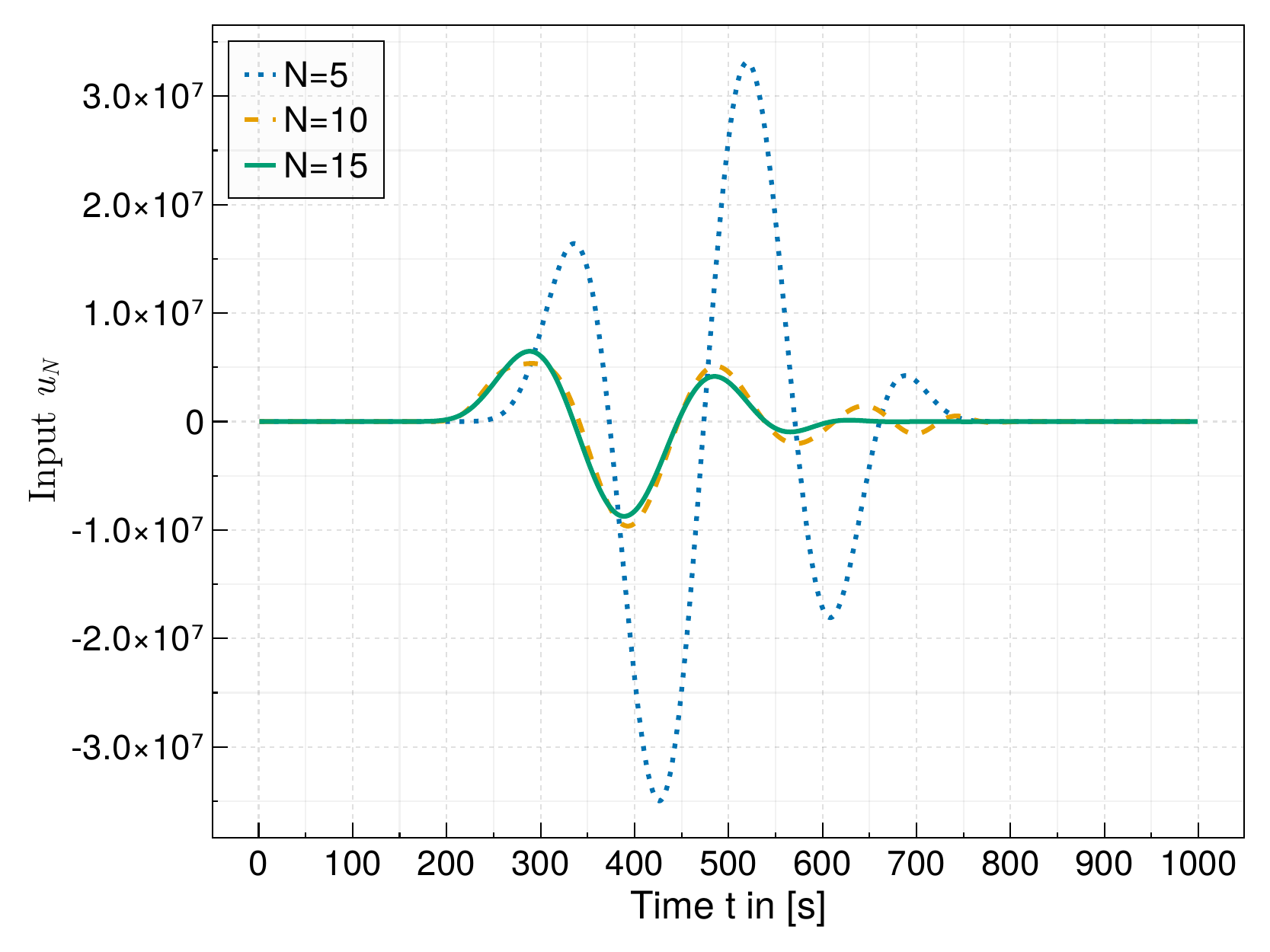}}
	\caption{Progress of approximated input signals $u_{N}$ for aluminum (left) and steel 38Si7 (right) with $T=1000$ and $\omega=2$.}
	\label{fig:input_signals_progress}
\end{figure}

We find an approximation of the signal input 
\begin{align}
	u(t) \approx&~ \frac{\lambda ~ \Delta y}{	\hat{\Omega}_{\omega,T}} ~ \sum_{i=0}^{N} \eta_{i} ~ \Omega_{\omega,T}^{(i)}(t) =: u_{N}(t) \label{eq:input_signal_approx}
\end{align}
where $N \in \mathbb{N}_{\geq 0}$ denotes the upper limit of iterations. In case of aluminum, we find a suitable approximation for $N=7$. In case of steel 38Si7, we need a higher number of iterations, e.g. $N=15$. The progress of input signals for aluminum with $N\in\{1,3,7\}$ and steel 38Si7 with $N\in\{5,10,15\}$ is presented in Figure \ref{fig:input_signals_progress}. These results also confirm our previous considerations: the input signal for steel 38Si7 shows more oscillatory behavior and significantly higher amplitudes than the input signal for aluminum.

In a nutshell, we find four important parameters which influence the computation of the input signal in Eq. \eqref{eq:fbc_input_signal}: length of rod $L$ and diffusivity $\alpha$, as part of sequence $\eta_{i}$, and final time $T$ and steepness $\omega$ in the derivatives of trajectory $\Phi_{\omega,T}$. 

\section{Simulation Results}
\label{sec:simulation_results}
In this section, we compare the computed input signals and the resulting heat conduction simulation for aluminum and steel 38Si7. We assume the physical properties as listed in Table \ref{table:physical_prop_al_steel} and the trajectory parameters $T=1000$ seconds and $\omega=2.0$.  We find $\hat{\Omega}_{\omega,T} \approx 17.06 \cdot 10^{-6}$ with Eq. \eqref{eq:omega_integral}, we assume the initial temperature $\heat_{0}(x) = 300$ and specify a temperature rise of $\Delta y = 100$ Kelvin. A maximum number of iterations $N=40$ is considered to approximate the input signal \eqref{eq:input_signal_approx} for both scenarios (aluminum and steel 38Si7) to yield a high-quality approximation of $N=\infty$. As explained in Section \ref{sec:influence_control}, a lower number of iterations than $N=40$ shall suffice for aluminum and steel 38Si7 as well. Heat equation \eqref{eq:fbc_heat_eq} is discretized in space using finite differences with $101$ grid points and is simulated using the integration method \textit{KenCarp4} (see \cite{article:kennedy2003additive}) from the \textsc{Julia} library \mbox{\textit{OrdinaryDiffEq.jl}} \cite{software:rackauckas2020sciml}. The plots are created with \textit{Makie.jl} \cite{article:danisch2021makie}.

\begin{figure}[!ht]
	\centering
	\subfloat[Aluminum]{\includegraphics[width=0.48\columnwidth]{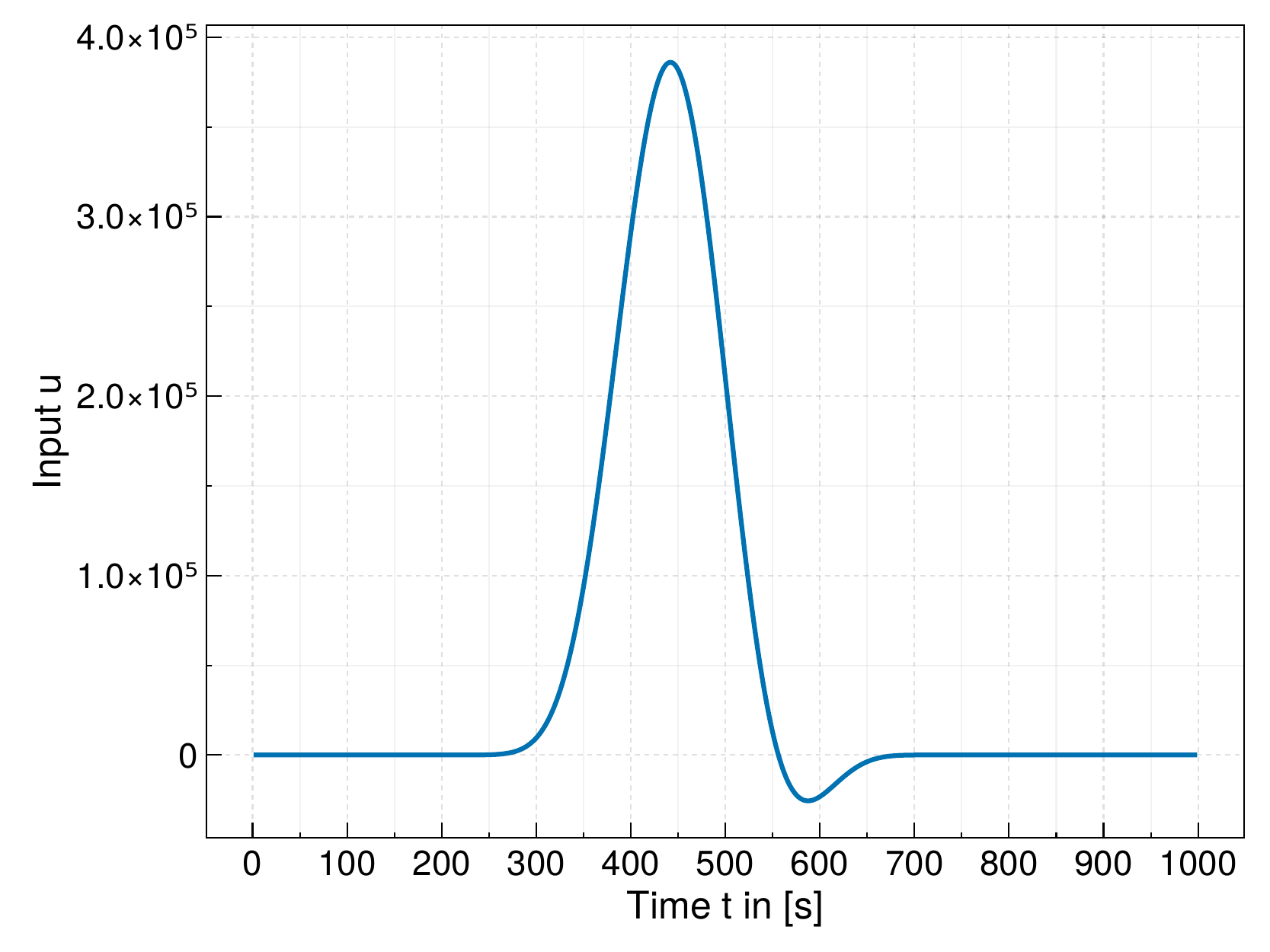}}
	\subfloat[Steel 38Si7]{\includegraphics[width=0.48\columnwidth]{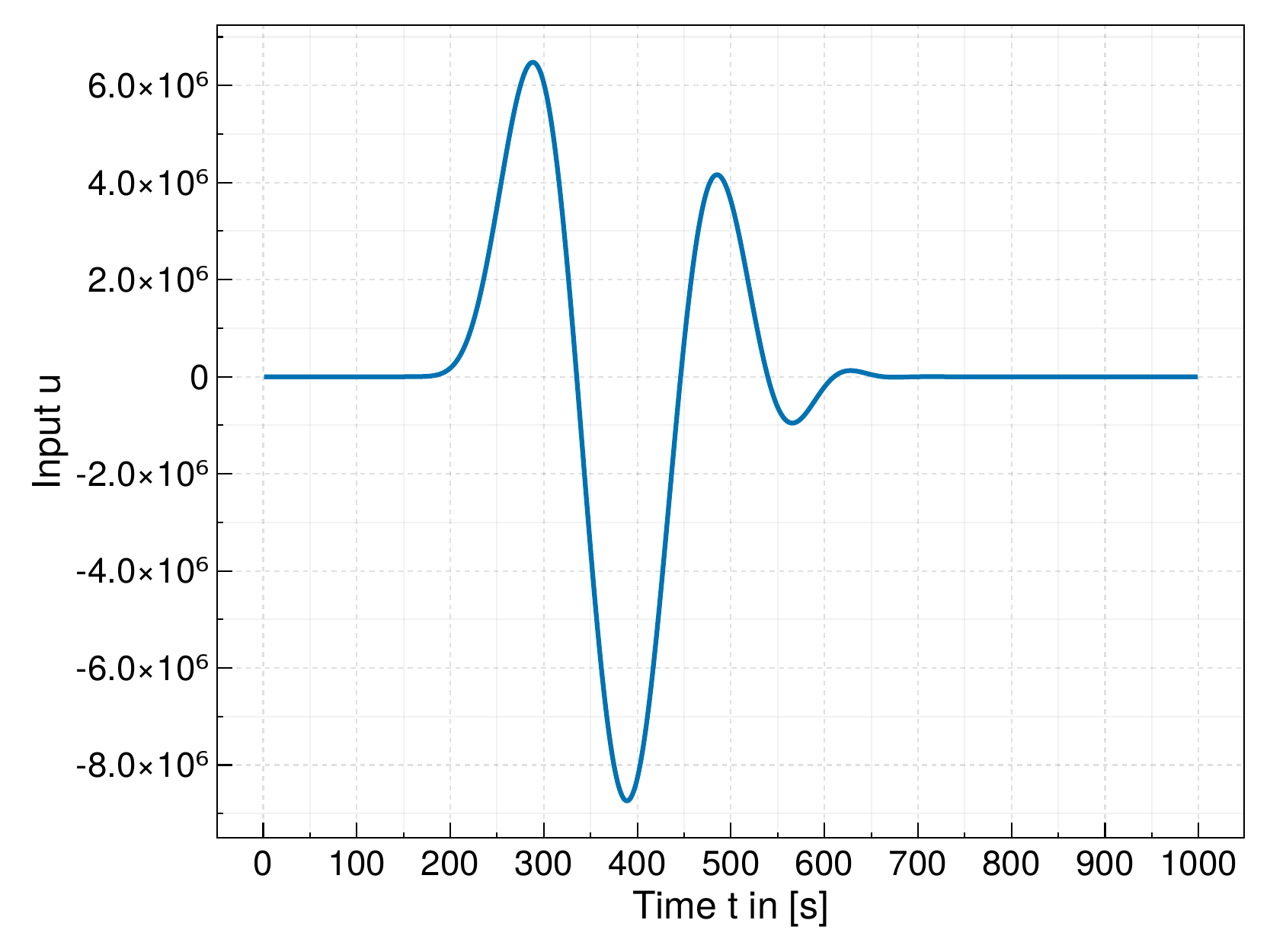}}
	
	\subfloat[Aluminum]{\includegraphics[width=0.48\columnwidth]{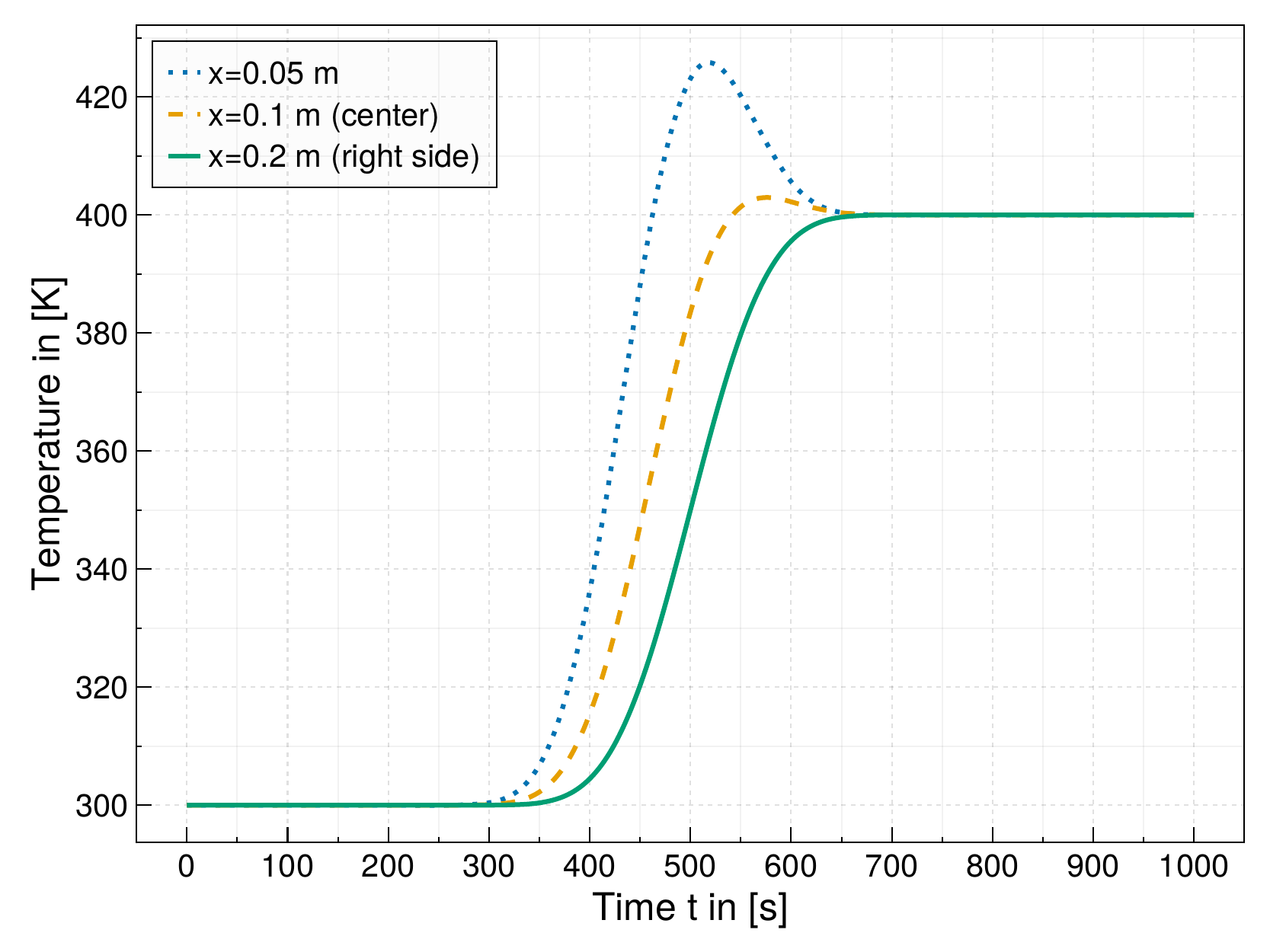}}
	\subfloat[Steel 38Si7]{\includegraphics[width=0.48\columnwidth]{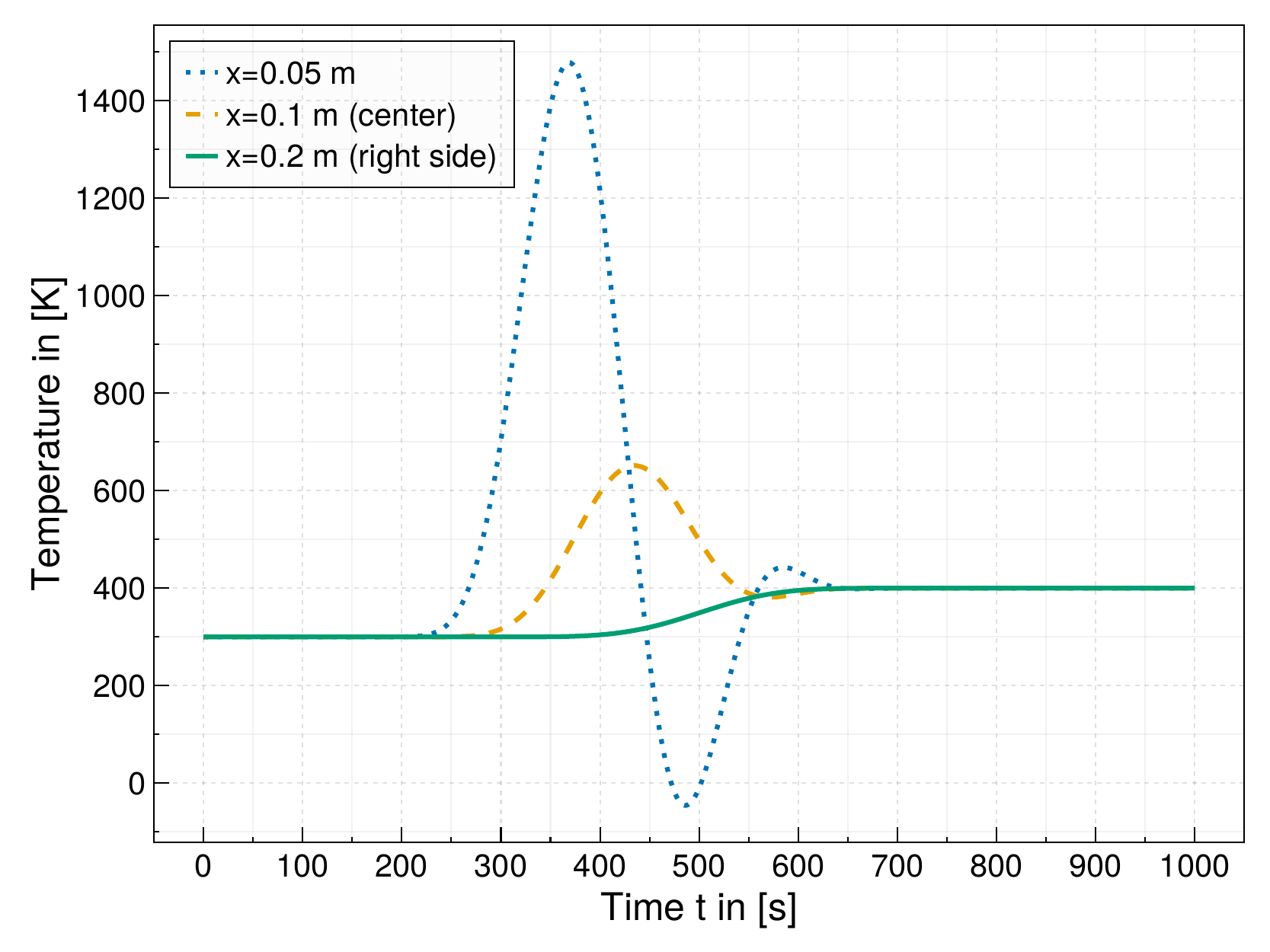}}
	\caption{Input signals and the resulting temperatures at position $x\in\{0.05,0.1,0.2\}$.}
	\label{fig:input_signals_temperatures}
\end{figure}

The input signals and the resulting temperatures are illustrated in Figure \ref{fig:input_signals_temperatures} for aluminum (a,c) and steel 38Si7 (b,d). In both cases the output signal, which is the measured temperature at $x=0.2$ meter, follows the reference and reaches $400$ Kelvin. Thus, from a pure \textit{mathematical} point of view, the input signals are computed correctly for both scenarios. However, from a \textit{physical} or \textit{technical} point of view, we need to discuss the input signals and the resulting temperatures rather critically. Firstly, it may not be possible to apply negative input signals if the actuator offers only heating and not cooling. Secondly, it is physically not possible to reach temperatures below zero Kelvin as portrayed for steel in Figure \ref{fig:input_signals_temperatures} (d). Therefore, the control parameters - final time $T$ and steepness $\omega$ - have to be readjusted as discussed in Section \ref{sec:influence_control} to yield physically sensible results. 

\section*{Source Code}
The source code is developed in \textsc{Julia} programming language and is available on \textit{GitHub}: 
\begin{center}
	\url{https://github.com/stephans3/BenchmarkFlatnessbasedControl.jl}
\end{center}

\section*{Conclusion}
In this article, we presented the computation of input signals for trajectory planning of a one-dimensional heat equation using flatness-based control design. We analyzed the influence of system and control parameters on the computation of the input signal. We found, that different material properties (aluminum, steel 38Si7) result in  completely different input signals and open-loop dynamics even if all other parameters (length of rod, final time, steepness of transition) are the same. Moreover, our calculations and discussions demonstrate the complexity of flatness-based control of simplified realistic heat conduction problems. Further research on the evaluation of flatness-based control design with focus on rather realistic scenarios in two and three dimensions is necessary to gain a deeper insight of this complexity.


\begin{thebibliography}{9}

	\bibitem{article:fliess1995flatness}
	Michel Fliess, Jean L\'evine, Philippe Martin, Pierre Rouchon:
	\emph{Flatness and defect of non-linear systems: introductory theory and examples}. International Journal of Control 61.6 (1995): 1327--1361.

	\bibitem{article:laroche2000motion}
	B\'eatrice Laroche, Philippe Martin, Pierre Rouchon: \emph{Motion planning for the heat equation}. International Journal of Robust and Nonlinear Control: IFAC‐Affiliated Journal 10.8 (2000): 629--643.
	
	\bibitem{article:ollivier2001generalization}
	Fran\c{c}ois Ollivier, Alexandre Sedoglavic: \emph{A generalization of flatness to nonlinear systems of partial differential equations. Application to the command of a flexible rod}. IFAC Proceedings Volumes 34.6 (2001): 219--223.
	
	\bibitem{book:rudolph2003flatness}
	Joachim Rudolph, Jan Winkler, Frank Woittennek:	\emph{Flatness Based Control of Distributed Parameter Systems: Examples and Computer Exercises from Various Technological Domains}. Shaker, 2003.
	
	\bibitem{article:utz2010trajectory}
	Tilman Utz, Knut Graichen, Andreas Kugi: \emph{Trajectory planning and receding horizon tracking control of a quasilinear diffusion-convection-reaction system}. IFAC Proceedings Volumes 43.14 (2010): 587--592.

	\bibitem{article:herzog2023an}
	Roland Herzog, Dmytro Strelnikov: \emph{An optimal control problem for single-spot pulsed laser welding}. Journal of Mathematics in Industry 13.1 (2023): 4.
	
	\bibitem{article:scholz2020modeling}	
	Stephan Scholz, Lothar Berger: \emph{Modeling of a multiple source heating plate}. arXiv preprint. arXiv:2011.14939 (2020).
	
	\bibitem{article:scholz2022hestia}	
	Stephan Scholz, L. Berger: \emph{Hestia.jl: A Julia library for heat conduction modeling with boundary actuation}. Simulation Notes Europe SNE 33.1 (2023): 27--30.

	\bibitem{book:baehr2013heat} 
	Hans Dieter Baehr, Karl Stephan: \emph{Heat and mass transfer}. Springer Science \& Business Media, 2013.
	
	\bibitem{book:lienhard2020heat} 
	John H. Lienhard IV, John H. Lienhard V: \emph{A Heat Transfer Textbook}. Phlogiston Press, 2020.
	
	\bibitem{online:periodic2023aluminium}
	Periodic Table: \emph{Aluminium – Periodic Table}. [Online]. Available:\url{https://www.periodic-table.org/Aluminium-periodic-table/}. [Accessed: Feb. 3, 2023].
	
	\bibitem{online:ovako202338Si7}
	Ovako: \emph{38Si7}. [Online]. Available: \url{https://steelnavigator.ovako.com/steel-grades/38si7/}. [Accessed: Feb. 3, 2023].
	
	
	\bibitem{software:fischer2021coni}
	Ferdinand Fischer, Jakob Gabriel, Simon Kerschbaum: \emph{coni-a Matlab toolbox facilitating the solution of control problems}. Zenodo (2021). Available: \url{https://zenodo.org/record/6420876}.

	\bibitem{software:scholz2023bellbruno} 
	Stephan Scholz: \emph{BellBruno.jl}. Zenodo. 2023. Available: \url{https://doi.org/10.5281/zenodo.7685927}


	\bibitem{article:kennedy2003additive}
	Christopher A. Kennedy, Mark H. Carpenter: \emph{Additive Runge-Kutta schemes for convection-diffusion-reaction equations}.
	Applied numerical mathematics 44.1-2 (2003): 139--181.

	\bibitem{software:rackauckas2020sciml} 
	Rackauckas C, contributors: \emph{SciML/DifferentialEquations.jl}. Zenodo. 2022. Available: \url{https://doi.org/10.5281/zenodo.7239171}
	
	
	\bibitem{article:danisch2021makie}
	Simon Danisch, Julius Krumbiegel: \emph{Makie.jl: Flexible high-performance data visualization for Julia}. Journal of Open Source Software 6.65 (2021): 3349.
	

\end{thebibliography}
\end{document}